\documentclass[notitlepage,11pt,reqno]{amsart}
\usepackage{amssymb,amsmath,amscd,xspace,tocvsec2,color,mathrsfs,enumerate}
\usepackage[breaklinks=true]{hyperref}
\usepackage[letterpaper]{geometry}
\theoremstyle{plain}

\newtheorem{theorem}[equation]{Theorem}

\theoremstyle{definition}
\newtheorem{defn}[equation]{Definition}
\newtheorem{remark}[equation]{Remark}

\numberwithin{equation}{section}

\newcommand{\wt}{\operatorname{wt}}
\newcommand{\conv}{\operatorname{conv}}

\newcommand{\ch}{\operatorname{ch}}

\newcommand{\vla}{\ensuremath{\mathbb{V}^\lambda}}

\newcommand{\R}{\ensuremath{\mathbb{R}}}
\newcommand{\C}{\ensuremath{\mathbb{C}}}

\newcommand{\bba}{\ensuremath{\mathbb{A}}}

\newcommand{\Z}{\ensuremath{\mathbb{Z}}}
\newcommand{\calp}{\ensuremath{\mathcal{P}}\xspace}
\newcommand{\liehr}{\ensuremath{\mathfrak{h}_{\mathbb{R}}}}

\newcommand{\lie}[1]{\ensuremath{\mathfrak{#1}\xspace}}
\newcommand{\comment}[1]{}

\begin{document}
\title[Weights of simple highest weight modules over a complex semisimple
Lie algebra]{Weights of simple highest weight modules\\ over a complex
semisimple Lie algebra}

\author[Apoorva Khare]{Apoorva Khare\\Stanford University}
\email[A.~Khare]{\tt khare@stanford.edu}
\address{Department of Mathematics, Stanford University, Stanford, CA -
94305}
\date{\today}
\thanks{The author was partially supported by the DARPA Grant YFA
N66001-11-1-4131.}


\begin{abstract}
In this short note we announce three formulas for the set of weights of
various classes of highest weight modules $\vla$ with highest weight
$\lambda$, over a complex semisimple Lie algebra $\lie{g}$ with Cartan
subalgebra $\lie{h}$. These include, but are not restricted to, all
(highest weight) simple modules $L(\lambda)$. We also assert that these
formulas are the ``best possible", in that they do not hold in general
for other highest weight modules in a very precise sense.

The proofs of the results in this note are included in an updated copy
(Version 3) of the paper\hfill\break
\href{http://arxiv.org/abs/1301.1140}{\tt
http://arxiv.org/abs/1301.1140}.
The proofs involve studying the convex hull of the set of
$\lie{h}$-weights $\wt \vla$ in their own right. Thus, we show that if
$\vla = L(\lambda)$ is simple, or if $\lambda$ is not on a simple root
hyperplane and $\vla$ is arbitrary, the hull of the infinite set $\wt
\vla$ is a convex polyhedron - i.e., cut out by only finitely many
hyperplanes.
(This extends the notion of the Weyl polytope to arbitrary simple modules
$L(\lambda)$.) It is also shown that the partially ordered set (under
quotienting) of modules $\vla$ with fixed convex hull, has unique
``largest" and ``smallest" elements.
\end{abstract}
\maketitle

In this note we study the structure of highest weight modules over a
complex semisimple Lie algebra $\lie{g}$. These modules, especially Verma
modules $M(\lambda)$ and their simple quotients $L(\lambda)$, have been
objects of much longstanding interest in representation theory.
Determining the structure (e.g., Jordan-Holder factors, the set of
weights, weight multiplicities) are important problems that have long
been studied in the field.

We provide three formulas to compute the set of weights of simple modules
$L(\lambda)$ for arbitrary weights $\lambda$. One of these formulas is in
terms of finite-dimensional submodules for a distinguished Levi
subalgebra, while another is in terms of the convex hull of the weights.
A consequence of these results is the extension of the Weyl Character
Formula to certain infinite-dimensional simple highest weight modules. It
is further shown that these formulas hold more generally for other
highest weight modules, but do not hold in general. These results (and
several other reasons from the literature) lead to our studying the
convex hull of the weights.

\section{Weights of (simple) highest weight modules}\label{S1}

Let $\lie{g}$ be a complex semisimple Lie algebra, with a fixed
triangular decomposition $\lie{g} = \lie{n}^- \oplus \lie{h} \oplus
\lie{n}^+$. Let $\Phi^\pm$ be the corresponding sets of positive and
negative roots, with simple roots and Weyl group given by $\Delta = \{
\alpha_i : i \in I \}$ and $W$ respectively.
Now suppose that $(,)$ is the nondegenerate bilinear form on $\lie{h}^*$
that is induced by the Killing form on $\lie{h}$, and $h_i \in \lie{h}$
corresponds via this form to $2 \alpha_i / (\alpha_i,\alpha_i)$ for all
$i \in I$. Also fix a set of Chevalley generators $\{ x_{\alpha_i}^\pm :
i \in I \}$ of $\lie{g}$ such that $[x^+_{\alpha_i}, x^-_{\alpha_j}] =
\delta_{i,j} h_i$.
Next, given a weight $\lambda \in \lie{h}^*$, let $M(\lambda) = U\lie{g}
/ U\lie{g}(\lie{n}^+ \oplus \ker_{\lie{h}} \lambda)$ denote the Verma
module with highest weight $\lambda$, and $L(\lambda)$ its unique simple
quotient. Given a $\lie{g}$-module $M$, let $\wt M$ denote the set of
$\lie{h}$-weights of $M$ (e.g., of quotients $\vla$ of $M(\lambda)$).
Finally, if $X \subset \lie{h}^*$ and $R \subset \R$, then let $R_+ := R
\cap [0,\infty)$, $\conv_\R X$ be the convex hull of $X$, and $RX$ denote
the set of finite linear combinations of $X$ with coefficients in
$R$.

In this note, we are interested in the structure of highest weight
modules $M(\lambda) \twoheadrightarrow \vla$. For instance, it is
well-known that for dominant integral highest weights $\lambda$ (i.e.,
$\lambda(h_i) \in \Z_+$ for all $i \in I$), the simple module
$L(\lambda)$ is finite-dimensional. Such modules and their weights yield
symmetries, combinatorial formulas, and crystals.
The convex hull of the set of weights is called the {\it Weyl polytope}
$\calp(\lambda) := \conv_\R \wt L(\lambda)$, and it is a bounded convex
$W$-invariant polyhedron, since $\wt L(\lambda)$ is a finite $W$-stable
set. Moreover (see \cite[Theorem 7.41]{Ha}),
\begin{equation}\label{Eweyl}
\wt L(\lambda) = (\lambda - \Z \Delta) \cap \conv_\R \wt L(\lambda).
\end{equation}

Another set of weights for which the structure of $\wt L(\lambda)$ is
easily computed is when $\lambda$ is antidominant - i.e., $(\alpha,
2\lambda + \sum_{\beta \in \Phi^+} \beta)/(\alpha,\alpha) -1 \notin \Z_+$
for every positive root $\alpha$. This is a Zariski dense set in
$\lie{h}^*$, being the complement of a countable set of affine
hyperplanes. In this case, $M(\lambda)$ is simple (see \cite[Theorem
4.8]{H3}), and hence is the unique highest weight module with highest
weight $\lambda$. Thus, $\wt L(\lambda) = \wt M(\lambda) = \lambda - \Z_+
\Delta$, and one checks that this equals $(\lambda - \Z \Delta) \cap
\conv_\R \wt L(\lambda)$.

These formulas naturally lead us to the main motivation for this note:
\begin{itemize}
\item To compute the set of weights of all simple modules $L(\lambda)$.

\item To determine whether Equation \eqref{Eweyl} holds more generally
for all $\lambda \in \lie{h}^*$.
\end{itemize}

Some progress towards these and related questions is known. For instance,
as per \cite[Chapter 7]{H3}, translation functors can be used to reduce
computing the formal character of simple modules $L(\lambda)$ for
semisimple $\lie{g}$ to the principal blocks $\mathcal{O}_0$ for all
$\lie{g}$. A more involved approach uses Kazhdan-Lusztig polynomials; see
\cite[Chapter 8]{H3} for a comprehensive treatment of this subject.

We now provide a complete resolution of both of the above questions, and
in greater generality. To do so, some additional notation is needed.

\begin{defn}
Fix $\lambda \in \lie{h}^*$ and $J \subset I$.
\begin{enumerate}
\item Define $J_\lambda := \{ i \in I : 0 \leq \lambda(h_i) \in \Z \}$.

\item Given a highest weight module $M(\lambda) \twoheadrightarrow \vla$
with the weight space $\vla_\lambda$ spanned by $v_\lambda$, define
\begin{equation}
J(\vla) := \{ i \in J_\lambda : (x^-_{\alpha_i})^{\lambda(h_i) + 1}
v_\lambda = 0 \}.
\end{equation}

\item Define $\Delta_J := \{ \alpha_i : i \in J \}$, and $W_J$ to be the
parabolic subgroup of $W$ generated by the simple
reflections corresponding to $\Delta_J$.

\item Denote the fundamental weights corresponding to the simple roots by
$\Omega := \{ \omega_i : i \in I \}$. Now define $\rho_J := \sum_{i \in
J} \omega_i$.

\item Define $\lie{g}_J$ to be the semisimple Lie subalgebra of
$\lie{g}$, generated by $\{ x^\pm_{\alpha_i} : i \in J \}$.

\item Given $J' \subset J_\lambda$, define the corresponding {\it
parabolic (or generalized) Verma module} to be
\[ M(\lambda,J') := U \lie{g} / U \lie{g} (\lie{n}^+ \oplus
\ker_{\lie{h}} \lambda \oplus \sum_{i \in J'} \C \cdot
(x^-_{\alpha_i})^{\lambda(h_i) + 1}). \]

\item Given $J' \subset J_\lambda$, define $L_{J'}(\lambda)$ to be the
finite-dimensional highest weight $\lie{g}_{J'}$-submodule of
$M(\lambda)$, which is generated by the highest weight space
$M(\lambda)_\lambda$. In other words, $L_{J'}(\lambda) = U \lie{g}_{J'}
\cdot M(\lambda)_\lambda$.
\end{enumerate}
\end{defn}

\begin{remark}
The set $J(\vla)$ is of fundamental importance in investigating the
behavior of $\vla$. Indeed, it can be shown using results in
\cite[Chapter 9]{H3} that $\vla$ is $\lie{g}_{J(\vla)}$-integrable, and
hence $\wt \vla$ is $W_{J(\vla)}$-stable. See \cite[Theorem 1]{Kh} for
more details. Additionally, $J(\vla)$ is closely related to the
classification theory for simple $\lie{h}$-weight $\lie{g}$-modules; see
\cite{Fe}.

Moreover, it is not hard to show (see again {\it loc.~cit.}) that if
$\vla$ is a Verma module $M(\lambda)$ or its simple quotient
$L(\lambda)$, then $J(\vla) = \emptyset, J_\lambda$ respectively. In the
special case that $\lambda$ is dominant integral and $\dim L(\lambda) <
\infty$, $J(L(\lambda)) = J_\lambda = I$.
Also note that parabolic Verma modules encompass all Verma modules (where
$J' = \emptyset$) and finite-dimensional simple modules (where $J' =
J_\lambda = I$). In fact, given $\lambda \in \lie{h}^*$ and $J' \subset
J_\lambda$, one has that $J(M(\lambda,J')) = J'$, and $L_{J'}(\lambda)$
is a simple finite-dimensional $\lie{g}_{J'}$-module.
\end{remark}

Here is the main result in this note. It provides explicit formulas to
compute the set of weights of all simple modules $L(\lambda)$, all
parabolic Verma modules $M(\lambda,J')$, and other modules $\vla$ as
well.

\begin{theorem}\label{Twtgvm}
Fix $\lambda \in \lie{h}^*$, and suppose $M(\lambda) \twoheadrightarrow
\vla$ is such that $|J_\lambda \setminus J(\vla)| \leq 1$, or $\vla =
M(\lambda,J')$ for some $J' \subset J_\lambda$. Then,
\begin{equation}\label{Ewtgvm}
\wt \vla = (\lambda - \Z \Delta) \cap \conv_\R \wt \vla = \wt
L_{J(\vla)}(\lambda) - \Z_+ (\Phi^+ \setminus \Phi^+_{J(\vla)}) =
\coprod_{\mu \in \Z_+ \Delta_{I \setminus J(\vla)}} \wt
L_{J(\vla)}(\lambda - \mu).
\end{equation}
\end{theorem}

\noindent Note that the last expression in Equation \eqref{Ewtgvm}
corresponds to the $\lie{g}_{J(\vla)}$-integrability of $\vla$ (as
discussed above), and it as well as the first equality both extend
\cite[Theorem 7.41]{Ha} to all simple modules and parabolic Verma
modules.

As an application of the techniques used in proving Theorem \ref{Twtgvm},
the Weyl Character Formula can be extended to various
infinite-dimensional simple modules $L(\lambda)$ as well. For instance,
the following result emerges out of our approach, although it was
(essentially) known in the literature - see \cite[Chapter 9]{H3}. To our
knowledge, the remaining theorems in this note are novel.

\begin{theorem}[Weyl Character Formula]\label{Twcf}
Suppose the set $S_\lambda := \{ w \bullet \lambda \ : \ w \in W, w
\bullet \lambda \leq \lambda \}$ equals $W_{J_\lambda}$. Then
$L(\lambda)$ is the unique quotient of $M(\lambda)$ whose set of weights
is $W_{J_\lambda}$-invariant, whence
\begin{equation}\label{Ewcf}
\ch L(\lambda) = \frac{\sum_{w \in W_{J_\lambda}} (-1)^{\ell(w)}
e^{w(\lambda + \rho_I)}}{\sum_{w \in W} (-1)^{\ell(w)} e^{w(\rho_I)}}.
\end{equation}
\end{theorem}

\noindent Just like Theorem \ref{Twtgvm}, this theorem also unifies the
cases of dominant integral and antidominant $\lambda$ (where $S_\lambda =
W$ and $J_\lambda = I$, or $S_\lambda = \{ 1 \}$ and $J_\lambda =
\emptyset$ respectively).

Note by Theorem \ref{Twtgvm} that the set of weights of an arbitrary
highest weight module $\vla$ can now be determined by computing the
multiplicities $[\vla : L(w \bullet \lambda)]$ of its Jordan-Holder
factors (which lie in the BGG Category $\mathcal{O}$).
A more direct attempt to compute $\wt \vla$ is to prove an analogue of
Theorem \ref{Twtgvm} for $\vla$. However, this result fails to hold for
all $\vla$. More precisely, Theorem \ref{Twtgvm} is the ``best possible"
when stated in terms of $|J_\lambda \setminus J(\vla)|$, in the following
sense.

\begin{theorem}\label{Twtsimple}
Fix $\lambda \in \lie{h}^*$, $M(\lambda) \twoheadrightarrow \vla$, and
$J' \subset J_\lambda$. If $|J_\lambda \setminus J'| \leq 1$, then
\begin{equation}\label{Ewtsimple}
J(\vla) = J' \quad \implies \quad \wt \vla = (\lambda - \Z \Delta) \cap
\conv_\R \wt \vla.
\end{equation}

\noindent Equation \eqref{Ewtsimple} however need not always hold if
$|J_\lambda \setminus J'| = 2$; and if $|J_\lambda \setminus J'| \geq 3$,
then Equation \eqref{Ewtsimple} always fails to hold for some $\vla$ with
$J(\vla) = J'$.
\end{theorem}

\noindent Thus, $\wt \vla$ does not always equal $\coprod_{\mu \in \Z_+
\Delta_{I \setminus J(\vla)}} \wt L_{J(\vla)}(\lambda - \mu)$
either. Another consequence of Theorem \ref{Twtsimple} is that the convex
hull $\conv_\R \wt \vla$ does not uniquely determine the module $\vla$
(or even its set of weights).

Theorems \ref{Twtgvm}, \ref{Twcf}, and \ref{Twtsimple} appear in an
updated copy (Version 3) of \cite{Kh} as Theorems 5, 6.5, and 6.2
respectively, and are proved therein. The techniques used in these proofs
involve representation theory and convexity theory.
Specifically, we study the convex hull of $\wt \vla$ and its properties,
for highest weight modules $M(\lambda) \twoheadrightarrow \vla$. In fact,
the above theorems - and various other motivations and results from the
literature - are strong reasons to also study convex hulls of weights.
It is not surprising that results for the convex hull are ``cleaner" than
those for the set of weights themselves.

\section{Convex hulls of weights of highest weight modules}

We now discuss some results concerning the convex hull of $\wt \vla$ for
several important families of highest weight modules $\vla$. First note
that all sets $\wt \vla$ are translates by $\lambda$ of subsets of the
root lattice $\Z \Delta$, so (via a translation) we always work in
Euclidean space $\liehr^* = \R \Delta$.
As mentioned above, if $\dim \vla < \infty$, then $\vla = L(\lambda)$ is
simple and $\lambda$ is dominant integral. In this case, $J(\vla) = I$
and $\calp(\lambda) = \conv_\R \wt \vla$ is a $W$-invariant compact
convex polyhedron. Moreover, the Weyl polytope $\calp(\lambda)$ has
vertex set $W(\lambda)$.
At the ``opposite end" are the Verma modules $\vla = M(\lambda)$ for all
$\lambda \in \lie{h}^*$, for which $J(M(\lambda)) = \emptyset$ and
$\conv_\R \wt M(\lambda)$ is a polyhedron with unique vertex $\lambda$.

Given these facts, it is natural to ask about the behavior of $\conv_\R
\wt \vla$ for various (infinite-dimensional) $\vla$. We show in a large
number of cases that the convex hull of $\wt \vla$ is in fact a convex
polyhedron, which is defined as a finite intersection of half-spaces
formed by affine hyperplanes. Note that it is not apparent that $\conv_\R
\wt \vla$ is cut out by only finitely many affine hyperplanes, especially
when $\dim \vla$ (and hence $\wt \vla$) is infinite.

The techniques used in proving the theorems stated in Section \ref{S1}
also yield many rewards in computing convex hulls. For instance,
\begin{itemize}
\item Computing the convex hulls of weights, for several families of
highest weight modules $\vla$, including all simple and parabolic Verma
modules. It is shown in all these cases that $\conv_\R \wt \vla$ is a
convex polyhedron.

\item Classifying the faces of these convex hulls, and writing down
precisely the inclusion relations between these faces.

\item Results in the literature (by Vinberg, Cellini, Chari, and others)
which were known earlier only for finite-dimensional simple modules, are
now shown for all highest weight modules.
\end{itemize}

Another feature of this work is that it focusses on several important
families of highest weight modules that feature prominently in the
literature:
\begin{enumerate}[(i)]
\item Parabolic Verma modules, which include all Verma modules as well as
finite-dimensional simple modules.
\item All simple highest weight modules $L(\lambda)$.
\item All modules $\vla$ with $\lambda$ not on a simple root hyperplane.
This is a larger Zariski dense set in $\lie{h}^*$ than the antidominant
weights $\lambda$ (for which $\vla = M(\lambda) = L(\lambda)$) as well as
all regular weights $\lambda$. We term these weights {\it
simply-regular}.
\end{enumerate}

\noindent We also consider a fourth class of highest weight modules
termed ``pure" modules. These modules feature in the classification of
all simple $\lie{h}$-weight $\lie{g}$-modules, in work of Fernando
\cite{Fe}. 

\begin{defn}\hfill
\begin{enumerate}
\item A $\lie{g}$-module $M$ is {\it pure} if for each $X \in \lie{g}$,
the set of $m \in M$ such that $\dim \C[X]m < \infty$ is either $0$ or
$M$.
\item Given a highest weight module $M(\lambda) \twoheadrightarrow \vla$
and $J \subset I$, define $\wt_J \vla := \wt \vla \cap (\lambda - \Z_+
\Delta)$.
\item (As above:) Define $\lambda \in \lie{h}^*$ to be {\it
simply-regular} if $(\lambda, \alpha_i) \neq 0$ for all $i \in I$.
\end{enumerate}
\end{defn}

The proofs of the results in Section \ref{S1} employ a wide variety of
techniques to study all of these families of modules. Thus, a module that
lies in more than one of these families can be studied in more than one
way. For instance, the following result holds for four different kinds of
highest weight modules. Corresponding to these, there are multiple
approaches and proofs given in \cite{Kh}.

\begin{theorem}\label{Tpolyhedron}
Suppose $(\lambda, \vla)$ satisfy one of the following:
(a) $\lambda \in \lie{h}^*$ is simply-regular and $\vla$ is arbitrary;
(b) $|J_\lambda \setminus J(\vla)| \leq 1$;
(c) $\vla = M(\lambda, J')$ for some $J' \subset J_\lambda$; or
(d) $\vla$ is pure.
\begin{enumerate}
\item The convex hull (in Euclidean space) $\conv_\R \wt \vla \subset
\lambda + \liehr^*$ is a convex polyhedron with vertices
$W_{J(\vla)}(\lambda)$.
\item The stabilizer subgroup in $W$ of both $\wt \vla$ and $\conv_\R \wt
\vla$ is $W_{J(\vla)}$.
\item The faces of $\conv_\R \wt \vla$ are precisely $\conv_\R w(\wt_J
\vla)$ for $J \subset I$ and $w \in W_{J(\vla)}$.
\item If $\lambda$ is simply-regular, the extremal rays at the vertex
$\lambda$ (i.e., the infinite length edges that pass through $\lambda$)
are $\{ \lambda - \R_+ \alpha_i \ : \ i \notin J(\vla) \}$.
\end{enumerate}
\end{theorem}

\noindent Thus, the results for $\conv_\R \wt \vla$ hold for a much
larger set of modules than the formula for $\wt \vla$ itself. It is also
easy to check that two faces are equal precisely when the set of weights
on them are equal. Thus we also prove a result that identifies more
generally for {\it all} highest weight modules $\vla$, the pairs $(w,J),
(w',J')$ for which $w(\wt_J \vla) = w'(\wt_{J'} \vla)$. See \cite[Theorem
4]{Kh}.

\begin{remark}
Theorem \ref{Tpolyhedron} holds for all simple modules $L(\lambda)$,
since $J(L(\lambda)) = J_\lambda$ and so (b) holds. Consequently, the
notion of the Weyl polytope extends to $L(\lambda)$ for arbitrary
$\lambda \in \lie{h}^*$, except that one now obtains a
$W_{J_\lambda}$-invariant convex polyhedron in general. Even more
generally, one can associate naturally to an arbitrary highest weight
module $\vla$, the $W_{J(\vla)}$-invariant polyhedron $\calp(\vla) :=
\conv_\R \wt M(\lambda, J(\vla))$.
\end{remark}

\begin{remark}
Theorem \ref{Tpolyhedron} was only known previously for parabolic Verma
modules - but that special case could not address the other ``nontrivial"
$\vla$, where $\lambda$ is not on countably many affine hyperplanes in
$\lie{h}^*$ (i.e., not antidominant). In light of Theorem
\ref{Tpolyhedron}, one can now work also with all simple modules
$L(\lambda)$, as well as all $\vla$ when $\lambda$ is not on the finite
set of simple root hyperplanes.
\end{remark}

To conclude this note, here is another application of our convexity- and
representation-theoretic approach. 
Note that the set of highest weight modules is naturally equipped with a
partial order under surjection, and it has unique maximal and minimal
elements $M(\lambda)$ and $L(\lambda)$ respectively. We now assert that
this ordering can be refined in terms of the stabilizer subgroup of the
weights, or equivalently, their convex hull.

\begin{theorem}\label{Tminmax}
Fix $\lambda \in \lie{h}^*$ and $J' \subset  J_\lambda$ such that either
$\lambda$ is simply-regular, or $J' = \emptyset$ or $J_\lambda$. There
exist unique ``largest" and ``smallest" highest weight modules
$M_{\max}(\lambda, J'), M_{\min}(\lambda,J')$ such that the following are
equivalent for a nonzero module $M(\lambda) \twoheadrightarrow \vla$:
\begin{enumerate}
\item $\conv_\R \wt \vla = \conv_\R \wt M(\lambda,J')$.
\item The stabilizer subgroup in $W$ of $\conv_\R \wt \vla$ is $W_{J'}$.
\item The largest parabolic subgroup of $W$ that preserves $\conv_\R \wt
\vla$ is $W_{J'}$.
\item $M(\lambda) \twoheadrightarrow M_{\max}(\lambda,J')
\twoheadrightarrow \vla \twoheadrightarrow M_{\min}(\lambda,J')$.
\end{enumerate}
\end{theorem}

\noindent As above, Theorem \ref{Tpolyhedron} appears as (parts of)
\cite[Theorems 2 and 3]{Kh}, and Theorem \ref{Tminmax} as \cite[Theorem
6]{Kh}. These results are also proved in {\it loc.~cit.}

\smallskip
\end{document}